\documentclass[AMA,STIX1COL]{WileyNJD-v2}

\newcommand\BibTeX{{\rmfamily B\kern-.05em \textsc{i\kern-.025em b}\kern-.08em
T\kern-.1667em\lower.7ex\hbox{E}\kern-.125emX}}

\articletype{Article Type}%

\received{<day> <Month>, <year>}
\revised{<day> <Month>, <year>}
\accepted{<day> <Month>, <year>}

\usepackage{mathtools}
\usepackage{bm}
\usepackage{nomencl}

\newcommand{\pp}[1]{\left(#1\right)}
\newcommand{\norme}[1]{\left|#1\right|}

\newcommand{\dd}[2]{\frac{\partial #1}{\partial #2}}

\newcommand{\Dt}{\Delta t}

{}

\usepackage{setspace}

\usepackage{graphicx}
\graphicspath{ {./im/} }

\begin{document}

\title{A High Order Stabilized Solver for the Volume Averaged Navier-Stokes Equations}

\author[1]{Toni El Geitani}

\author[1]{Shahab Golshan}

\author[1]{Bruno Blais*}

\authormark{Toni El Geitani \textsc{et al}}

\address[1]{\orgdiv{Research Unit for Industrial Flows Processes (URPEI), Department of Chemical
Engineering}, \orgname{École Polytechique de Montréal}, \orgaddress{\state{Quebec}, \country{Canada}}}

\corres{*Bruno Blais. \email{bruno.blais@polymtl.ca}}

\presentaddress{PO Box 6079, Stn Centre-Ville, Montréal,
QC, Canada, H3C 3A7}

\abstract[Abstract]{The Volume-Averaged Navier-Stokes equations are used to study fluid flow in the presence of fixed or moving solids such as packed or fluidized beds. We develop a high-order finite element solver using both forms A and B of these equations. We introduce tailored stabilization techniques to prevent oscillations in regions of sharp gradients, to relax the Ladyzhenskaya-Babuska-Brezzi inf-sup condition, and to enhance the local mass conservation and the robustness of the formulation. We calculate the void fraction using the Particle Centroid Method. Using different drag models, we calculate the drag force exerted by the solids on the fluid. We implement the method of manufactured solution to verify our solver. We demonstrate that the model preserves the order of convergence of the underlying finite element discretization. Finally, we simulate gas flow through a randomly packed bed and study the pressure drop and mass conservation properties to validate our model.}


\jnlcitation{\cname{%
\author{El Geitani T.}, 
\author{Golshan S.}, and
\author{Blais B}} (\cyear{2022}), 
\ctitle{A High Order Stabilized Solver for the Volume Averaged Navier-Stokes Equations}, \cjournal{Int J Numer Meth Fluids}, \cvol{}.}

\maketitle

\section{Introduction}\label{introduction}
Multi-phase flows are systems with two or more distinct phases flowing simultaneously within a mixture. Their classification depends on the combination of phases involved: solid-gas, solid-liquid, immiscible liquid-liquid, or liquid-gas flows \cite{ch1}. Solid-gas flows play a key role in the dynamics of fluidized and packed bed reactors, spouted beds, pneumatic conveyors and other technological applications. For the solid-fluid flows in most chemical engineering processes, it is well known that the phenomena that occur at the particle scale greatly influence the macroscopic behavior of the process \cite{vanderHoefM.A2006MMoG}. Generally, the design and scale-up of these beds are based on expensive and time-consuming pilot-scale tests. Therefore, modeling has become an important topic for researchers due to its contribution to our understanding of the processes involved as well as its predictive capabilities of solid-fluid flows in engineering-scale equipment. It has numerous advantages such as allowing the understanding of various phenomena occurring in processes, enabling sensitivity analysis on different input parameters, and testing several configurations and operational conditions at a much cheaper cost than experimental work \cite{norouzi2016,GolshanShahab2017Hoss}.
In order to study such systems, the volume-averaged Navier Stokes (VANS) equations are used as they ensure that the volume of fluid that can occupy a cell and the volume of solids occupying it respect the continuity equation, thus ensuring mass conservation. Since in CFD-DEM, the fluid is modeled at a scale larger than that of the particles (mesoscopic), the governing equations for the fluid phase (VANS) are obtained based on local averaging methods as in the Two Fluid Model (TFM) \cite{GolshanShahab2020Raio}.

Mainly, there exist two forms of the VANS equations: Model A (Set II) and Model B (Set I) \cite{GIDASPOW19941, ZHOUZ.Y2010Dpso}. Model A is derived from Model B by decomposing the particle-fluid interactions into two components: a component accounting for the macroscopic variations in the fluid stress tensor on a large scale, and another component accounting for the variations in the point stress tensor when the fluid passes around particles \cite{ZHOUZ.Y2010Dpso}. It is the most widely used set of equations as it is implemented in several commercial software such as FLUENT and CFX \cite{ZHOUZ.Y2010Dpso}. Gidaspow derived a Model C which replaces the solid velocity in Model B with the relative velocity. According to Gidaspow \cite{GIDASPOW19941}, in order to be invariant under a change of frame of reference, entropy production should be a function of relative velocity and not of individual phase velocities. Zhou et al. \cite{ZHOUZ.Y2010Dpso} simplified Model B (Set I) under the assumption $\rho_f(1-\epsilon_f)[\dd{\bm{u}_f}{t} + \nabla \cdot(\bm{u}_f\bm{u}_f)]=0$ in order to obtain Set III. The main difference between these models is that the solid-fluid interaction forces in Set I are explicit while in Set III, they are implicit. However, Set III is not general but conditional as it is only valid when the assumption is satisfied \cite{norouzi2016}. This assumption does not hold for numerous solid-fluid flows of interest and, as such, this model should not be considered. 

Many numerical methods exist for solving these equations, but the majority of solvers for the VANS equation use the Finite Volume Method (FVM) \cite{finitevolume} or the Lattice Boltzmann Method (LBM) \cite{lbm, BlaisBruno2015AclB}. However, the Finite Element Method (FEM) is an attractive candidate to solve these equations because of its inherent capacity for higher order formulations. It is considered a generalization of the classical variational (the Ritz) and weighted-residual (such as weak-form Galerkin Least Square (GLS), etc.) methods.


 High-order formulations are methods possessing an order of accuracy greater than two \cite{WangZ.J2007Hmft}. When the error in the solution $e$ obtained with a method decreases with the size of the mesh $h$ such that $e \propto h^{k+1}$, then the method's error is said to be of order $k+1$. The use of low order methods is usually prevalent in simulations due to them being robust. Even though high-order methods can achieve higher accuracy for the same computational cost, they are less used due to their complexity and because they are often considered less robust \cite{HuynhH.T2014Hmfc}. Steady-state convergence is more difficult and challenging to reach with high-order simulation, and for worse cases, the solution may even diverge \cite{WangZ.J2007Hmft}. However, high-order methods have shown promise in efficiently solving problems with a requirement for high accuracy such as vortex dominated flows, numerical simulation of turbulence, and large-eddy simulations \cite{WangZ.J2007Hmft, LiZhaorui2012Ahfd}.
The advantage of these methods is their relatively small error \cite{highOrderAccuracy}. Moreover, for the same number of degrees of freedom (DOFs), high-order methods generally produce solutions with less errors compared to lower order methods. However, when dealing with explicit time-stepping methods, high-order schemes usually have a smaller critical CFL value than low order methods \cite{WangZ.J2007Hmft}. This critical CFL value is a condition for the stability of the underlying scheme. This issue is often overcome by using an implicit scheme.

In this paper, we build a high-order FEM stabilized VANS solver as a first step in realizing a coupled CFD-DEM solver within the open-source software Lethe \cite{BlaisBruno2020LAop}. We present a fully implicit finite element VANS solver. Thus, it does not require any stability condition. Moreover, the formulation we introduce is well posed in the presence of smooth void fraction or by using void fraction smoothing and it is globally mass conservative as it accounts for the time variation of the void fraction in the continuity equation. The order of our solver is arbitrary and we prove that it can support a large scale of high-order methods. 
In our solver, we calculate the void fraction in each cell using the Particle Centroid Method (PCM), then we project it onto the FEM mesh. When desired, we apply a smoothing factor to resolve the oscillations arising from the discontinuous void fraction and a primal-dual active set method to bound the void fraction. We use the continuity equation as the basis of the grad-div stabilization in order to improve local mass conservation. Using two different drag models, we calculate the drag force on individual particles which contributes to the particle-fluid interactions in the VANS equations. Through the method of manufactured solutions (MMS), we verify the solver by analyzing the order of convergence in space and time. Finally, we simulate a packed bed and study the pressure drop across the bed. We compare the results with the pressure drop obtained from the Ergun equation in order to validate our model.

\section{Governing Equations and Numerical Implementation}\label{implementation}
As forms A and B of the VANS equations are the most commonly known models, we represent them in Eqs. (\ref{eq::vans_1}), (\ref{eq::vans_21}) and (\ref{eq::vans_2}). The conservative form of the continuity equation for an incompressible fluid remains the same for both models:
\begin{eqnarray}
    \rho_f\dd{\pp{\epsilon_f}}{t} + \rho_f \nabla \cdot \pp{\epsilon_f \bm{u}} = \mathit{m'}
    \label{eq::vans_1}
\end{eqnarray}
where $\rho_f$ the fluid density, which is assumed to be constant, $\bm{u}$ is the fluid velocity, $\epsilon_f$ the fluid's void fraction, and $\mathit{m'}$ the volumetric source of mass \cite{GIDASPOW19941}.
The difference between the two models lies in the formulation of the momentum equations. The conservative form of Model A of the incompressible VANS equations is defined as:
\begin{eqnarray}
\rho_f \pp{\dd{\pp{\epsilon_f \bm{u}}}{t} + \nabla \cdot \pp{\epsilon_f \bm{u} \otimes \bm{u}}} = -\epsilon_f \nabla p + \epsilon_f \nabla \cdot \pp{\bm{\tau}_f} + \bm{F_{pf}} + \rho_f\epsilon_f \bm{f} 
\label{eq::vans_21}
\end{eqnarray}
The conservative form of Model B of the incompressible VANS equations is:
\begin{eqnarray}
\rho_f \pp{\dd{\pp{\epsilon_f \bm{u}}}{t} + \nabla \cdot \pp{\epsilon_f \bm{u} \otimes \bm{u}}} = -\nabla p + \nabla \cdot \pp{\bm{\tau}_f} + \bm{F_{pf}} + \rho_f\epsilon_f \bm{f} 
\label{eq::vans_2}
\end{eqnarray}
with 
\begin{eqnarray}
\bm{\tau}_f = \nu \pp{\pp{\nabla \bm{u}} + \pp{\nabla \bm{u}}^T} - \frac{2}{3} \nu \nabla \cdot \bm{u} \bm{I}\label{eq::stress}
\end{eqnarray}

where $p$ is  the pressure, $\bm{\tau}_f$ the deviatoric stress tensor, $\nu = \frac{\mu_f}{\rho}$ the kinematic viscosity, $\bm{f}$ the external force(ex. gravity), and $\bm{F}_{pf}$ is the momentum transfer term between the solid and fluid phases and includes forces such as drag, virtual mass, Basset force, Saffman lift, and Magnus lift \cite{BlaisBruno2015Otuo}. $\bm{I}$ is the identity or unit matrix. In the case of incompressible fluid, $\nabla \cdot \bm{\tau} = \nu \pp{\nabla^2 \bm{u}}$ in Eq. (\ref{eq::stress}).

Each form leads to a different pressure equation. In Form A, $\epsilon_f$ is taken into account on the continuous equation and alters the saddle point problem for the pressure since $\epsilon_f$ multiplies the pressure gradient directly while in form B, $\epsilon_f$ is lumped with the momentum coupling term \cite{BerardAriane2020Emic}. Both models are mathematically equivalent and differ in the way they implement incompressibility. Form B of the equations achieves a realistic representation of the system where the total pressure drop occurs in the fluid phase and no pressure drop occurs in the solid phase. This is not true for the other model \cite{GIDASPOW19941}. This is because Form B uses pressure at previous times to calculate subsequent pressure values. Also, Form B possesses all real characteristics rendering the model well-posed as described by Gidaspow \cite{GIDASPOW19941}. On the other hand, and according to Bouillard et al. \cite{BouillardJ.X1989Pdia}, Form A of the VANS equations represents an ill-posed initial value problem that possesses complex characteristics which make this model conditionally stable. We solved both Form A and Form B of the VANS equations and compared the results obtained from both models. 

In order to obtain the weak form of the VANS equations, we determine its non-conservative form. Using the incompressibility constraint, we simplify our problem even further as the density becomes constant.

\begin{eqnarray}
    \dd{\bm{\epsilon_f}}{t} + \nabla \cdot  \epsilon_f \bm{u} = \frac{\mathit{m'}}{\rho_f} \label{eq::vans_1_non_conservative}
\end{eqnarray}
Form B: 
\begin{eqnarray}
\label{eq::vans_2_non_conservativeB}
& \epsilon_f  \dd{\bm{u}}{t} +  \epsilon_f \bm{u} \cdot \nabla \bm{u} + \bm{u} \Bigg(\dd{\bm{\ \epsilon_f}}{t} + \nabla \cdot \epsilon_f \bm{u} \Bigg) = 
\\& \frac{1}{\rho_f} \Big(-\nabla p + \nabla \cdot \pp{\bm{\tau_f}} + \bm{F_{pf}} \Big) + \epsilon_f \bm{f} \notag
\end{eqnarray}
Form A:
\begin{eqnarray}
\label{eq::vans_2_non_conservativeA}
& \epsilon_f  \dd{\bm{u}}{t} +  \epsilon_f \bm{u} \cdot \nabla \bm{u} + \bm{u} \Bigg(\dd{\bm{\ \epsilon_f}}{t} + \nabla \cdot \epsilon_f \bm{u} \Bigg) = 
\\& \frac{1}{\rho_f} \Big(-\epsilon_f\nabla p + \epsilon_f\nabla \cdot \pp{\bm{\tau_f}} + \bm{F_{pf}} \Big) + \epsilon_f \bm{f} \notag
\end{eqnarray}

Since the continuity equation is used along with the momentum equation to solve the fluid flow, the left hand side term in brackets of Eqs. (\ref{eq::vans_2_non_conservativeB}) and (\ref{eq::vans_2_non_conservativeA}), which represents the continuity equation can be replaced by the mass source term.  We also expand the term $\pp{\nabla \cdot \epsilon_f\bm{u}}$.

\begin{eqnarray}
    \dd{\bm{\epsilon_f}}{t} + \epsilon_f\nabla \cdot\bm{u} + \bm{u}\nabla \epsilon_f = \frac{\mathit{m'}}{\rho_f} \label{eq::vans_1_non_conservative2}
\end{eqnarray}
Form B: 
\begin{eqnarray}
\label{eq::vans_2_non_conservative_2}
&\epsilon_f  \dd{\bm{u}}{t} + \epsilon_f \bm{u} \cdot \nabla \bm{u} + \bm{u} \Bigg( \underbrace{\dd{\bm{\epsilon_f}}{t} + \epsilon_f\nabla \cdot\bm{u} + \bm{u}\nabla \epsilon_f}_\mathit{m'}\Bigg) = 
\\& \frac{1}{\rho_f}\pp{-\nabla p + \nabla \cdot \pp{\bm{\tau_f}} + \bm{F_{pf}}} + \epsilon_f \bm{f} \notag
\end{eqnarray}
Form A:
\begin{eqnarray}
\label{eq::vans_2_non_conservative_2A}
&\epsilon_f  \dd{\bm{u}}{t} + \epsilon_f \bm{u} \cdot \nabla \bm{u} + \bm{u} \Bigg( \underbrace{\dd{\bm{\epsilon_f}}{t} + \epsilon_f\nabla \cdot\bm{u} + \bm{u}\nabla \epsilon_f}_\mathit{m'}\Bigg) = 
\\& \frac{1}{\rho_f}\pp{-\epsilon_f\nabla p + \epsilon_f\nabla \cdot \pp{\bm{\tau_f}} + \bm{F_{pf}}} + \epsilon_f \bm{f} \notag
\end{eqnarray}

The equations were solved by the finite element method using stabilization techniques which allows to circumvent the limitations of the classical Galerkin approach and satisfies the Ladyzhenskaya–Babuska–Brezzi (LBB) inf–sup condition. For a more detailed explanation, we refer the reader to the Lethe paper \cite{BlaisBruno2020LAop, Tezduyar2003}.

\subsection{The Weak Form and Finite Element Equations}
We implement a stabilized formulation in order to solve the VANS equations. Two stabilizing terms, streamline-upwind/Petrove-Galerkin (SUPG) and pressure-stabilizing/Petrov-Galerkin (PSPG) were added. PSPG helps relax the LBB condition and allows the use of equal order elements while the SUPG term stabilizes unresolved boundary layers \cite{BlaisBruno2020LAop}. The grad-div stabilization was added to ensure local mass conservation. 

\subsubsection{The Stabilized Formulation}
The detailed weak form of the VANS equations can be expressed by setting $\nabla \cdot \bm{\tau} = \nu \pp{\nabla^2 \bm{u}}$.
\begin{eqnarray}
&\int_{\Omega}\pp{\dd{\pp{\epsilon_f}}{t} + \epsilon_f\nabla \cdot\bm{u} + \bm{u}\nabla \epsilon_f} q d\Omega
+ \mathcal{SR} \cdot \pp{\tau_{u} \nabla q} d\Omega_k=0 \label{eq::VANS_gls1}
\end{eqnarray}
Form B:
\begin{eqnarray}
&\int_{\Omega} \pp{\epsilon_f  \dd{\bm{u}}{t} + \epsilon_f \bm{u} \cdot \nabla \bm{u}} \cdot \bm{v} d\Omega + \frac{1}{\rho_f} \Bigg(\int_{\Omega} \nu \pp{\nabla \bm{u}}\nabla \bm{v} d\Omega +  \int_{\Omega} \mathit{m'}\bm{u} \cdot \bm{v}d\Omega 
\\&- \int_{\Omega} p \nabla \cdot \bm{v} d\Omega +\int_{\Omega}\frac{\bm{F_{pf}}}{V_{\Omega}} \cdot \bm{v} d\Omega \Bigg) \nonumber  - \int_{\Omega}\epsilon_f \bm{f} \cdot \bm{v} d\Omega + \sum_K \int_{\Omega_k} \mathcal{SR} \cdot  \pp{\tau_{u}\bm{u}\cdot \nabla 
\bm{v}} d\Omega_k = 0 \label{eq::VANS_gls2B}
\end{eqnarray}
and
\begin{eqnarray}
    \mathcal{SR} =  \epsilon_f \dd{\bm{u}}{t} + \epsilon_f \bm{u} \nabla \cdot \bm{u} + \frac{1}{\rho_f}\Bigg(\mathit{m'}\bm{u}  + \nabla p - \nu \pp{\nabla^2 \bm{u}}  + \frac{\bm{F_{pf}}}{V_{\Omega}}\Bigg) - \epsilon_f \bm{f}
\end{eqnarray}
Form A:
\begin{eqnarray}
&\int_{\Omega} \pp{\epsilon_f  \dd{\bm{u}}{t} + \epsilon_f \bm{u} \cdot \nabla \bm{u}} \cdot \bm{v} d\Omega + \frac{1}{\rho_f} \Bigg(\int_{\Omega} \Big(\epsilon_f \nu \pp{\nabla^2 \bm{u}}+ \nu \nabla \bm{u}\nabla\epsilon_f \Big)\cdot \bm{v} d\Omega +  \int_{\Omega} \mathit{m'}\bm{u} \cdot \bm{v}d\Omega 
\\&- \int_{\Omega} \Big(\epsilon_f p \cdot \nabla \bm{v} + p \nabla\epsilon_f \cdot \bm{v} \Big) d\Omega +\int_{\Omega}\frac{\bm{F_{pf}}}{V_{\Omega}} \cdot \bm{v} d\Omega \Bigg) \nonumber  - \int_{\Omega}\epsilon_f \bm{f} \cdot \bm{v} d\Omega + \sum_K \int_{\Omega_k} \mathcal{SR} \cdot  \pp{\tau_{u}\bm{u}\cdot \nabla 
\bm{v}} d\Omega_k = 0 \label{eq::VANS_gls2A}
\end{eqnarray}
and
\begin{eqnarray}
    \mathcal{SR} =  \epsilon_f \dd{\bm{u}}{t} + \epsilon_f \bm{u} \nabla \cdot \bm{u} + \frac{1}{\rho_f}\Bigg(\mathit{m'}\bm{u}  + \epsilon_f\nabla p - \epsilon_f\nu \pp{\nabla^2 \bm{u}}   + \frac{\bm{F_{pf}}}{V_{\Omega}}\Bigg) - \epsilon_f \bm{f}
\end{eqnarray}
and $\mathcal{SR}$ is the strong residual associated with the PSPG and SUPG stabilizations. $\bm{q}$ and $\bm{v}$ are the test functions of the finite element method. $V_{\Omega}$ is the volume of the domain $\Omega$. The stabilization parameter $\tau_{u}$ is defined by Tezduyar \cite{Tezduyar2003} and is explained in detail in the Lethe paper \cite{BlaisBruno2020LAop}. Two forms of $\tau_{u}$ were implemented depending on whether the equations solved are stationary or time dependent. For transient simulations, $\tau_{u}$ is defined as:
\begin{eqnarray}
\tau_u = \left[ \pp{\frac{1}{\Dt}}^2 + \pp{\frac{2\norme{\bm{u}}}{h_{conv}}}^2 + 9\pp{\frac{4\nu}{h_{diff}^2}}^2 \right]^{-1/2}
\end{eqnarray}
while for stationary problems, $\tau_{u}$ is:
\begin{eqnarray}
\tau_u = \left[ \pp{\frac{2\norme{\bm{u}}}{h_{conv}}}^2 + 9\pp{\frac{4\nu}{
h_{diff}^2}}^2 \right]^{-1/2}
\end{eqnarray}
where $\Dt$ is the time step, $h_{conv}$ and $h_{diff}$ are the size of the element related to the convection transport and diffusion mechanism respectively \cite{BlaisBruno2020LAop}.

The weak form of Model A has been derived in another work \cite{Srivastava}. However, the derivation is incomplete as the weak form of the pressure formulation lacks a term. The weak form of the pressure term is defined as:
\begin{eqnarray}
    \int_{\Omega} \Big(\epsilon_f\nabla p \Big)\cdot \bm{v} d\Omega = \int_{\Omega} \Big(\epsilon_f p \cdot \nabla \bm{v} + p \nabla\epsilon_f \cdot \bm{v}\Big)d \Omega
\end{eqnarray}
The second term on the right hand side is missing in the weak form of Srivastava et al. \cite{Srivastava}. Even though this term does not appear to have any physical meaning, its presence is important to make the equation mathematically correct and equivalent to that of Model B as well as to obtain the correct pressure drop in coupled simulations.

\subsubsection{Grad-div Stabilization}
In order to improve mass conservation in incompressible flow simulations and to stabilize our system at higher $Re$ values, we apply grad-div stabilization to the VANS equations similarly to how it was applied to the Navier Stokes equations in Lethe \cite{BlaisBruno2020LAop}. However, the void fraction should be accounted for when considering the VANS equations. The resulting stabilization term added to the weak form of the VANS momentum equation is expressed as:
\begin{eqnarray}
   \sum_K \int_{\Omega_k} \gamma \pp{\dd{\epsilon_f}{t} + \nabla \cdot \pp{\epsilon_f \bm{u}}} \pp{\nabla \cdot\bm{v}} d\Omega_k
\end{eqnarray}
where $\gamma$ is a parameter related to the augmented Lagrangian formulation \cite{BlaisBruno2020LAop}.  This $\gamma$ parameter is defined differently by multiple authors in the literature. We chose to implement the formulation derived by Olshanskii et al. \cite{OLSHANSKII20093975} where $\gamma = \nu + \bm{u}h$ $m^2/s$ where h is the size of the element. The system we are solving is generally stiff and higher $Re$ values make it more difficult to solve. A constant value of $\gamma$ was initially used where we noticed that it didn't completely help in solving the issue, especially at a high $Re$ number. Therefore, the formulation given by  by Olshanskii et al. \cite{OLSHANSKII20093975} is used to stabilize the fluid flow. The appendix shows a simple test case with a sharp void fraction gradient which results in shocks in the velocity field. This was achieved by using a step function for the void fraction. We compare the velocity fields obtained using different combinations of stabilization terms. Before we implemented the grad-div stabilization, mass was not correctly conserved locally at the location of void fraction discontinuity where we had large velocity gradients. The application of this stabilization technique helps enhance the local conservation of mass. This conservation becomes much better as we refine the simulation mesh. This can be attributed to the element dependent $\gamma$. Refining the mesh also results in better oscillation control and reduces the distance over which the oscillations propagate.

\subsection{Void Fraction Calculation} 
We consider two cases for the calculation of the void fraction. In the case of the manufactured solution (MMS), we calculate the void fraction from an analytical expression that satisfies the VANS equation. In the case where particles exist, we calculate the void fraction in a finite element as:
\begin{eqnarray}
\epsilon_f = \frac{v_{element} - v_{particles}}{v_{element}}
\label{eq::void_fraction}
\end{eqnarray}
where $v_{element}$ is the volume of the finite element and $v_{particles}$ is the volume of all particles present within the element. In order to determine the volume of particles in an element, we use the particle centroid method (PCM) even though it possesses several limitations, due to its simplicity \cite{GolshanShahab2020Raio}.

\subsubsection{The Particle Centroid Method}
The particle-centered or particle-centroid method assumes the volume of the particle to be located at its centroid. However, this method entails some disadvantages mainly the large fluctuation of the void fraction with time when particles' centers cross cell boundaries \cite{ClarkeDanielA2018IoVF}. 
Another important limitation we should account for is the ratio of cell size to particle size which should be large enough in order to make the PCM accurate. But using a cell size which is much larger than particle size will prevent us from capturing the fluid flow field \cite{ClarkeDanielA2018IoVF}. Therefore, there is a trade-off between cell size and particle size in unresolved CFD-DEM. This latter issue can be mitigated by increasing the order of the scheme instead of refining the mesh.

\subsubsection{$\mathcal{L}^2$ Projection of Void Fraction} 
The void fraction obtained using the PCM method is calculated at each cell of the FEM mesh. However, for the assembly of the finite element system, we need the values of the void fraction at the nodes of the mesh. This is equivalent to finding the coefficients $\epsilon_j$ as to minimize the following equation \cite{larson}:
\begin{eqnarray}
    \min_{\epsilon \in \mathbb{R}} \frac{1}{2} \sum_i \Big(\sum_j\epsilon_j \varphi_j - f_i\Big)
    \label{L2_projection}
\end{eqnarray}
where $\epsilon_j$ is the projected void fraction on the finite element mesh and $f_i$ is the calculated void fraction in each cell.
We, therefore, applied $\mathcal{L}^2$ projection in order to project the calculated void fraction onto the quadrature points. This is done by assembling and solving the following finite element system:
\begin{eqnarray}
 \int_{\Omega} \varphi_i \epsilon_j\varphi_j d\Omega
  = \int_{\Omega} f_i \varphi_i  d\Omega
  \label{eq:l2}
\end{eqnarray}

The projection of the void fraction solves one of the limitations of the PCM method, where now the void fraction becomes continuous in space between elements.

However, the resulting values of the projected void fraction may exceed the value of 1 at certain locations (where there are strong gradients). Therefore, in order to ensure a realistic representation of the void fraction throughout the simulation domain, we bound the $\mathcal{L}^2$ system with a lower bound of 0.36, which represents the void fraction in a non-uniform packing and an upper bound of 1.

\subsubsection{Bounding of the Void Fraction}
The $\mathcal{L}^2$ projection is bounded using a primal-dual active set method detailed in step-41 of the deal.II tutorials where we introduce a Lagrange multiplier $\lambda$ in order to bound the problem \cite{dealII93}. 
This method allows us to keep the nature of the $\mathcal{L}^2$ projection as a quadratic optimization problem with n variables, but now with $n$ bound constraints on these variables. This is applicable in the case of bi/tri-linear elements. In the case of bi/tri-quadratic elements, this method becomes expensive and exhaustive since high-order elements do not have their minimum and maximum attained at the nodes.
We used the following two bounds for our problem:
\begin{eqnarray}
    \epsilon_{min} \leq \epsilon_j \leq \epsilon_{max}
\end{eqnarray}
where $\epsilon_{max} = 1$ is the maximum void fraction in the absence of particles while $\epsilon_{min}$ is the minimum void fraction and can be set by the user. This method, however, might reduce the mass conservation of our system as it slightly changes the mass preservative properties of the $\mathcal{L}^2$ projection. 

\subsubsection{Smoothing of the Void Fraction}
In certain applications, the void fraction can be discontinuous or semi-discontinuous. For example, in a packed bed, the void fraction at the entry and exit of the bed drops suddenly and rapidly from a value of 1 to a value of below 0.4 at the inlet of the bed and vice versa at the exit of the bed. This leads to a very high gradient in the void fraction which results in high velocity gradients and leads to undesired oscillations in the flow field at high Reynolds number. In order to circumvent this problem, we apply to the left hand side of Eq. (\ref{eq:l2}) a smoothing term similar to a Poisson equation to smooth the void fraction field. 
This smoothing term can be expressed as: 
\begin{eqnarray}
    \iint_{\Omega} L^2 \nabla \varphi_i \nabla \varphi_j d\Omega
\end{eqnarray}
where $L$ is the smoothing length chosen to be a multiple of the particle's diameter. In our case, $L$ is chosen to be twice the particle's diameter. The resulting void fraction distribution is sufficiently smooth resulting in an enhanced convergence of the non-linear solver. However, implementing grad-div stabilization allows us to reduce the smoothing length required.

\subsection{Fluid-Particle Interactions} 
As an initial implementation of our FEM-VANS code, we determine that the drag force is the only non-negligible force. This is because the Basset and virtual mass forces are only significant in unsteady flows while the lift force is not significant in packed beds \cite{GolshanShahab2020Raio}. However, since we are only simulating steady flow with a large number of particles, these forces become insignificant and can therefore be neglected. 

When many particles exist in the fluid, the influence of particles on the drag of a neighboring particle should be taken into account. The drag for this type of flow can be calculated using:
\begin{eqnarray}
    \bm{F_{D}} = \beta(\bm{u}_f-\bm{u}_p)
\end{eqnarray}
where $\bm{u}_p$ is the particle velocity, and $\beta$ is the interphase momentum exchange coefficient and is determined by several correlations found in the literature. Among these correlations, widely used models include the Gidaspow model \cite{GIDASPOW19941} which combines the Ergun equation\cite{Ergun} for $\epsilon_{f} < 0.8$ and the Wen-Yu model \cite{Wen} for $\epsilon_{f} \geq 0.8$, Syamlal and O'Brien model \cite{SyamlalM1988Sogl}, Di Felice Model \cite{DiFeliceR1994Tvff}, Koch and Hill model \cite{KochDonaldL2001IEIS}, Rong et al. model \cite{RongL.W2013Lsof} and Beestra et al. model \cite{BeetstraR2007Nsos}. For more information about the different drag models, we direct the reader to the article by Norouzi et al. \cite{NorouziHamidReza2021Otdf} and the article by Bérard et al.\cite{BerardAriane2020Emic}. Gidaspow explains that the interphase momentum exchange coefficient $\beta$ differs between models A and B according to \cite{GIDASPOW19941}:
\begin{eqnarray}
    \beta_{B} = \frac{\beta_A}{\epsilon_f}
    \label{eq::beta}
\end{eqnarray}
In our study, the drag force is calculated as a summation over particles of individual particle forces. We studied several drag models in order to evaluate the drag force. The drag force for Model B is calculated as:
\begin{eqnarray}
    \bm{F_{{pf}_B}} = \sum_{n=1}^{np} \frac{1}{2}\rho_f C_D A_{ref} \norme{\bm{u}_{f,p}-\bm{u}_p}(\bm{u}_{f,q}-\bm{u}_{p,avg})
\end{eqnarray}
and by using Eq. \ref{eq::beta}, we deduce that the drag force for model A is determined as :
\begin{eqnarray}
    \bm{F_{{pf}_B}} = \frac{\bm{F_{{pf}_A}}}{\epsilon_f}
\end{eqnarray}

where $n_p$ is the total number of particles in the cell, $\bm{u}_{f,p}$ is the interpolated fluid velocity at the particle's location calculated at the previous time step, $\bm{u}_{f,q} = \bm{u}_f$ is the fluid velocity at the quadrature point, and $\bm{u}_{p,avg}$ is the average particles' velocity in the cell both calculated at the current time step. This decomposition of the velocities allows us to consider the influence of the drag in the Jacobian of the matrix and makes the solid-fluid coupling semi-implicit thus resulting in faster and more robust convergence. $A_{ref}$ is the particle's reference area which is taken as the cross-sectional area:
\begin{eqnarray}
    A_{ref} = \pi r_p^2
\end{eqnarray}
where $r_p$ is the particle's radius and  $C_D$ is the drag coefficient and is determined depending on the model used. 

For the Di Felice model \cite{DiFeliceR1994Tvff}, the drag coefficient is:
\begin{eqnarray}
    C_{D_A} = \Big(0.63 + \frac{4.8}{\sqrt{Re_p}}\Big)^2 \epsilon_f^{2-\big[3.7 - 0.65 e^{\big(\frac{-(1.5-log_{10}(Re_p))^2}{2}\big)}\big]}
\end{eqnarray}
and $Re_p$ is the particle Reynolds number and is expressed as:
\begin{eqnarray}
    Re_p =\frac{\rho_f\norme{\bm{u}_f-\bm{u}_p}d_p}{\mu_f}
\end{eqnarray}
while for the Rong et al. model \cite{RongL.W2013Lsof}, the drag coefficient is defined as:
\begin{eqnarray}
    C_{D_A} = \Big(0.63 + \frac{4.8}{\sqrt{Re_p}}\Big)^2 \epsilon_f^{2-\big[2.65(\epsilon_f + 1) - (5.3-3.5\epsilon_f)\epsilon_f^2 e^{\big(\frac{-(1.5-log_{10}(Re_p))^2}{2}\big)}\big]}
\end{eqnarray}

For a packed bed, the particles are at rest and therefore, $\bm{u}_p = 0$.

\section{The Method of Manufactured Solutions for the Verification of the Stabilized VANS Scheme}\label{verification}
The Method of Manufactured Solutions (MMS) is a procedure which allows for the verification of a code's accuracy by manufacturing an analytical problem with an exact solution. Since the verification is purely mathematical, the solution does not have to be physically realistic; however, it has to be sufficiently complex in order to ensure a complete verification \cite{RoachePatrickJ2002CVbt}. With the help of a Python script and through symbolic manipulation, we created several test cases with increasing complexity in order to fully verify all aspects of the VANS equations. We adopted a similar approach as presented by Blais et al. \cite{BlaisBruno2015Otuo}. From Eqs. (\ref{eq::vans_1}), (\ref{eq::vans_21}), and (\ref{eq::vans_2}), we identify two source terms that should be accounted for in the MMS. The mass source term ($\mathit{m'}$) and the momentum source term $\bm{G}$. We define the momentum source term for each form as:
\begin{align}
\bm{G}_{B} = \rho_f \pp{\dd{\pp{\epsilon_f \bm{u}}}{t} + \nabla \cdot \pp{\epsilon_f \bm{u} \otimes \bm{u}}} - \mathit{m'}\bm{u} + \nabla p - \nabla \cdot \pp{\bm{\tau_f}} - \bm{F_{pf}} - \epsilon_f \bm{f} 
\label{eq::source_term}
\end{align}
\begin{align}
\bm{G}_{A} = \rho_f \pp{\dd{\pp{\epsilon_f \bm{u}}}{t} + \nabla \cdot \pp{\epsilon_f \bm{u} \otimes \bm{u}}} - \mathit{m'}\bm{u} + \epsilon_f\nabla p - \epsilon_f\nabla \cdot \pp{\bm{\tau_f}} - \bm{F_{pf}} - \epsilon_f \bm{f} 
\label{eq::source_termA}
\end{align}
No particles are considered in the MMS cases, and therefore, $\bm{F_{pf}} = 0$. Four test cases were created for the verification of our code. All test cases presented in the following section can be found in the lethe repository on Github at \url{https://github.com/lethe-cfd/lethe-utils/tree/master/cases/vans}.

We performed the simulations with a mesh ranging from 256 cells to 4096 cell on a fixed domain of simulation $\omega = [-1,1]\times[-1,1]$. Different orders for the velocity and pressure shape functions are used in order to determine whether or not our solver conserves the order of accuracy. Convergence is reached when the residual of the non-linear solver becomes inferior to $10^{-8}$. We calculate the order of convergence by a least-square linear regression of the Euclidean norm of the error with respect to the mesh size as was done by Blais et al. \cite{BlaisBruno2015Otuo}.

\begin{align}
    \norme{E}_2 = \sqrt{\int_\Omega \norme{e_q}^2}
\end{align}
where $\norme{E}_2$ is the $\mathcal{L}^2$ norm of the error and $e$ is the error at the quadrature points $q$.

We create four cases with increasing complexity to consider different possible scenarios. The $0^{th}$ case includes a constant void fraction and a steady state divergence-free velocity profile. Due to its simplicity, it is not presented in this paper. For all other cases, the pressure, void fraction and velocity fields are chosen to be sufficiently differentiable. The $1^{st}$ case investigates a spatially varying void fraction and a divergence free steady velocity. This is the only case that is studied for both Model A and Model B of the VANS equations to ensure both models are mathematically correct. Since the difference between the models is only related to the derivation of the pressure and stress tensor weak forms and in order to avoid repetition, we only present the other test cases for Model B. The $2^{nd}$ case investigates a non-divergence free steady velocity. The $3^{rd}$ case investigates a spatially and time varying void fraction as well as a transient velocity field.

The code was run for various finite element orders for velocity and pressure to determine whether or not the chosen order is well preserved within our implementation. The orders used were Q1-Q1, Q2-Q1, Q2-Q2, Q3-Q2, and Q3-Q3. The void fraction order is chosen to be equal to that of the velocity. For the transient case, and due to the large simulation time associated with higher orders, the study is limited to a highest order of Q2-Q2. 

\subsection{Case 1: Steady State Divergence-free Flow Problem}
The velocity ($\bm{u}$), pressure ($p$) and void fraction ($\epsilon$) for this case are defined as:
\begin{align}
\bm{u} =
\begin{bmatrix}
-2sin^2(\pi x)sin(\pi y)cos(\pi y) \\
2sin(\pi x)sin^2(\pi y)cos(\pi x) \\
0
\end{bmatrix}
\end{align}

\begin{align}
    p = sin(\pi x)sin(\pi y)
\end{align}

\begin{align}
    \epsilon = \frac{1}{2} + \frac{1}{4}sin(\pi x)sin(\pi y)
\end{align}

\begin{figure}[ht]
     \centering
         \includegraphics[width=0.7\textwidth]{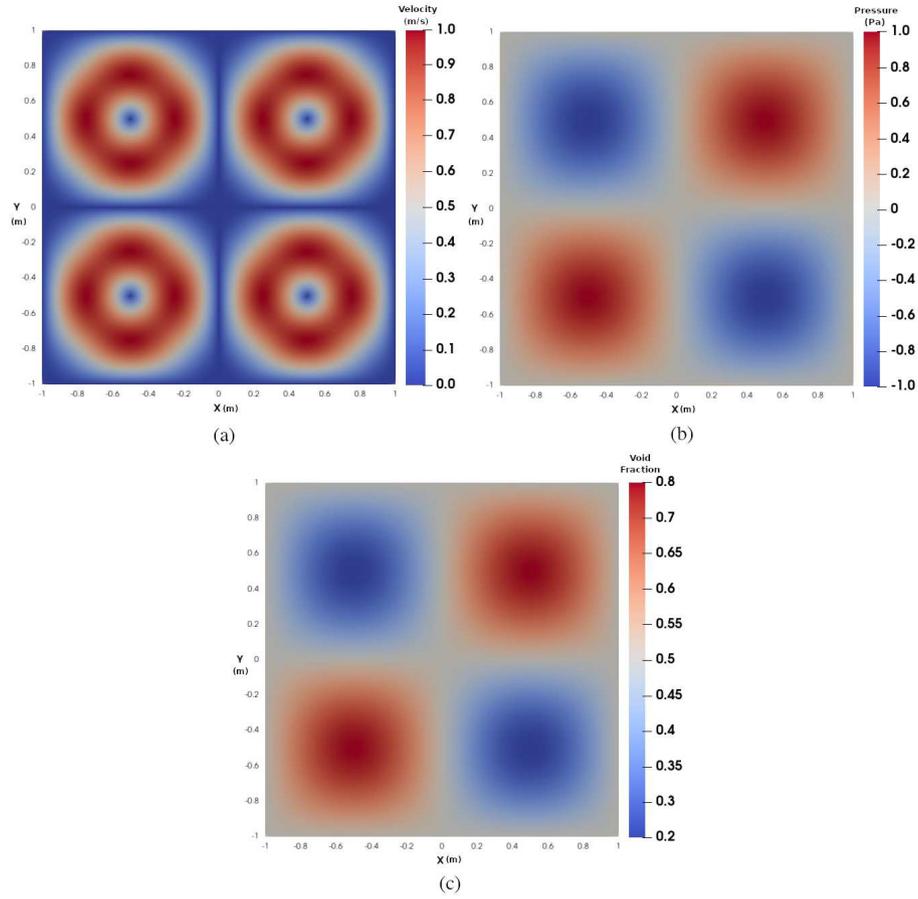}
        \caption{Case 1 analytical solution for (a) velocity, (b) pressure and (c) void fraction  using model B with 4096 cells and Q1-Q1 elements.}
        \label{fig:case1}
\end{figure}

\begin{figure}[ht]
     \centering
         \centering
         \includegraphics[width=0.7\textwidth]{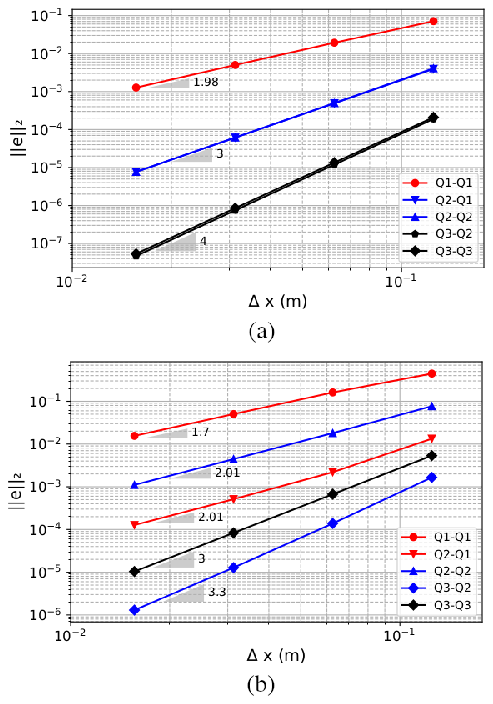}
        \caption{Case 1 convergence plots for (a) velocity and (b) pressure using Model B with respect to mesh size.}
        \label{fig:case1_convergence}
\end{figure}

\begin{figure}[ht]
     \centering
         \includegraphics[width=0.7\textwidth]{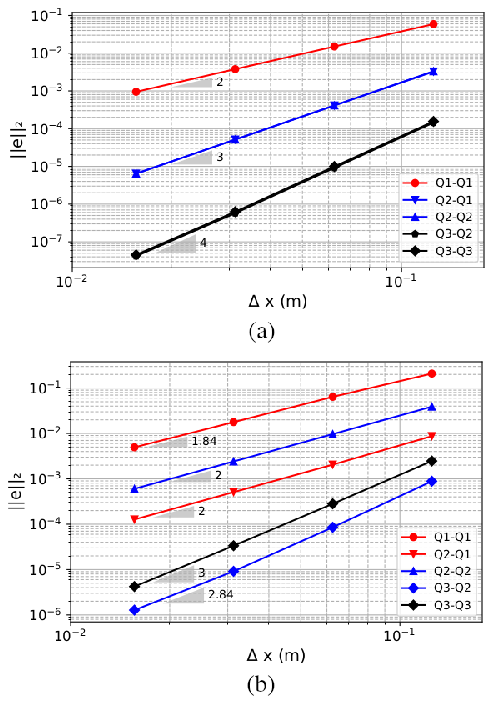}
        \caption{Case 1 convergence plots for (a) velocity and (b) pressure using Model A with respect to mesh size.}
        \label{fig:case1_convergenceA}
\end{figure}

In this case, there is no mass source term since mass conservation is satisfied. 
All fields are steady in this case, and the velocity field is divergence free which results in the normal stresses in the viscous stress tensor being zero. The velocity, pressure and void fraction fields are shown in Fig. \ref{fig:case1} for Model B. The convergence plots for models B and A are shown in Fig. \ref{fig:case1_convergence} and Fig. \ref{fig:case1_convergenceA} respectively.

\subsection{Case 2: Steady State Non Divergence-free Flow Problem}
The velocity, pressure and void fraction for this cases are defined as:
\begin{align}
\bm{u} =
\begin{bmatrix}
\frac{1}{e}e^{sin(\pi x)sin(\pi y)} \\
\frac{1}{e}e^{sin(\pi x)sin(\pi y)} \\
0
\end{bmatrix}
\end{align}

\begin{align}
    p = \frac{1}{2} + \frac{1}{2}sin(\pi x)sin(\pi y)
\end{align}

\begin{align}
    \epsilon = \frac{1}{e}e^{-sin(\pi x)sin(\pi y)}
\end{align}

In this case, mass conservation is inherently satisfied by the velocity field. Therefore, there is no mass source term. 

\begin{figure}[ht]
     \centering
         \includegraphics[width=0.7\textwidth]{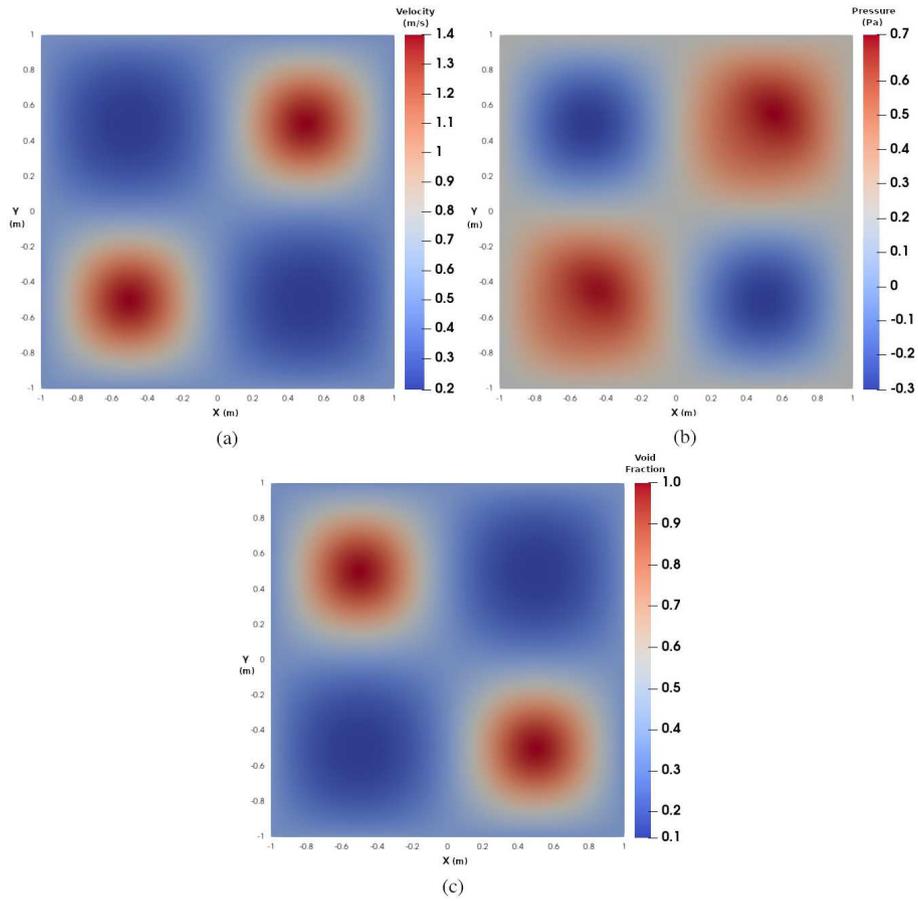}
        \caption{Case 2 analytical solution for (a) velocity, (b) pressure and (c) void fraction  using model B with 4096 cells and Q1-Q1 elements.}
        \label{fig:case2}
\end{figure}

\begin{figure}[ht]
     \centering
         \centering
         \includegraphics[width=0.7\textwidth]{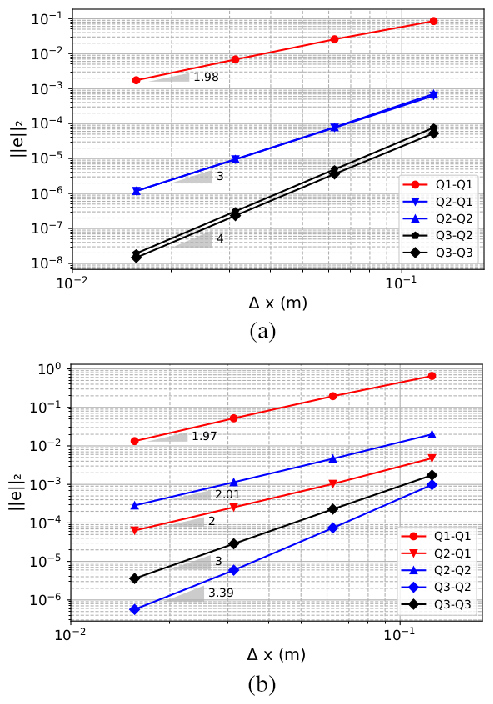}
        \caption{Case 2 convergence plots for (a) velocity and (b) pressure using Model B with respect to mesh size.}
        \label{fig:case2_convergence}
\end{figure}

All fields are steady in this case; however, the velocity field is not divergence free which results in the presence of all components of the stress tensor in the VANS equations. The velocity, pressure and void fraction fields are shown in Fig. \ref{fig:case2}. The convergence plot is shown in Fig. \ref{fig:case2_convergence}.

\subsection{Case 3: Unsteady Non Divergence-free Flow Problem}
This case is defined as:
\begin{align}
\bm{u} = cos(2\pi t)
\begin{bmatrix}
cos(\pi x)cos(\pi y) \\
cos(\pi x)cos(\pi y) \\
0
\end{bmatrix}
\end{align}

\begin{align}
    p = 0
\end{align}

\begin{align}
    \epsilon = \frac{1 - 0.1cos(2\pi t)e^{-sin(\pi x)sin(\pi y)}}{e}
\end{align}

\begin{figure}[ht]
     \centering
         \includegraphics[width=0.7\textwidth]{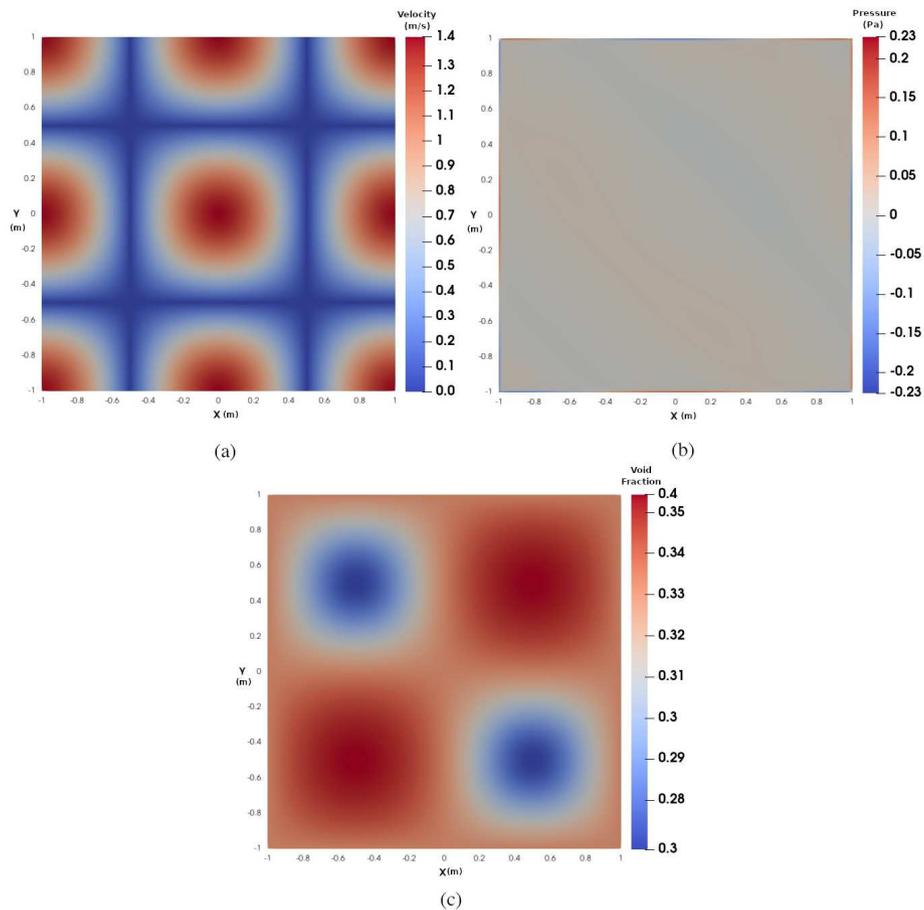}
        \caption{Case 3 analytical solution for (a) velocity, (b) pressure and (c) void fraction  using model B with 4096 cells and Q1-Q1 elements.}
        \label{fig:case3}
\end{figure}

\begin{figure}[ht]
     \centering
         \centering
         \includegraphics[width=0.7\textwidth]{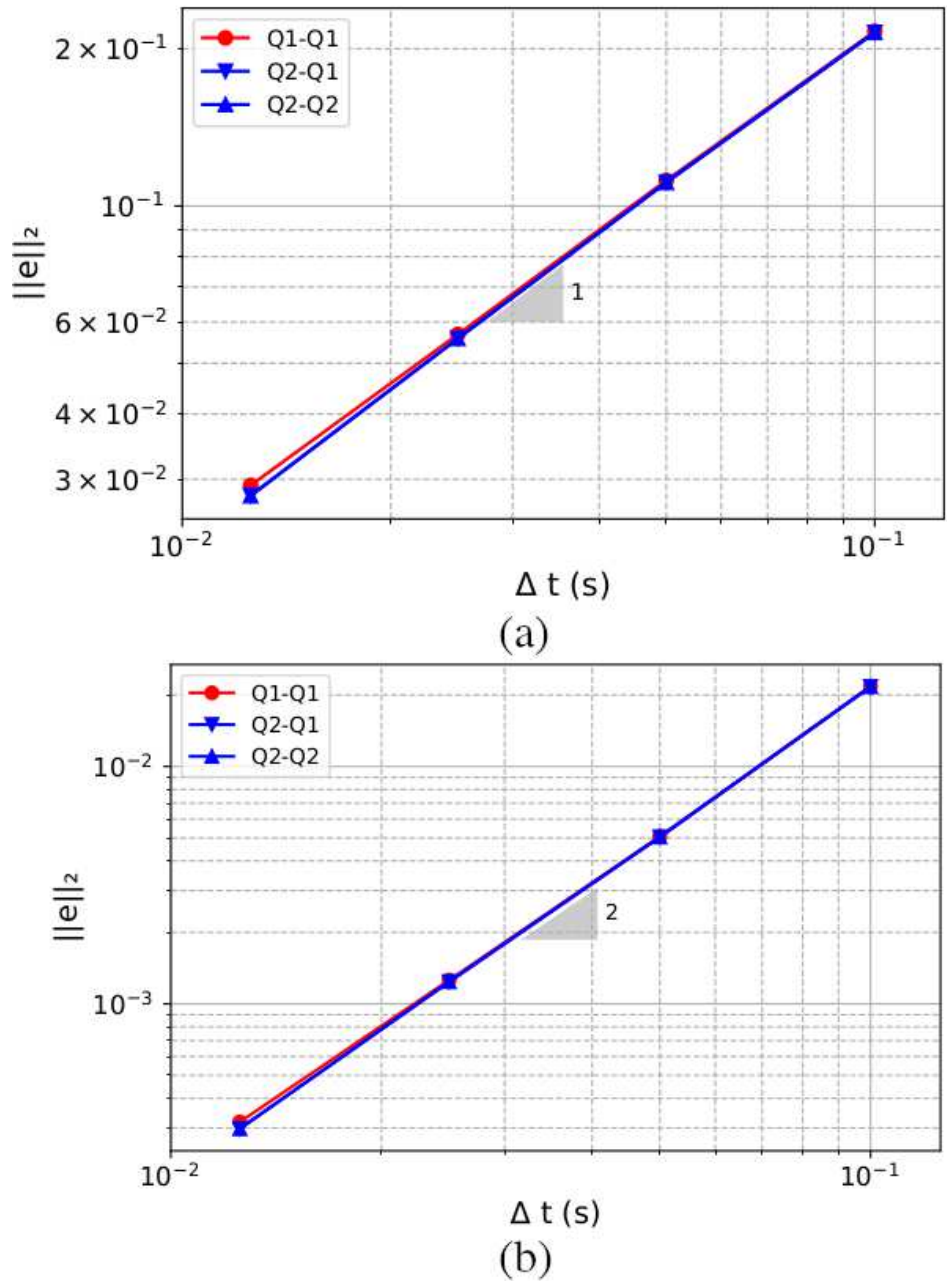}
        \caption{Case 3 velocity convergence plots for (a) BDF1 and (b) BDF2 using Model B with respect to mesh size.}
        \label{fig:case3_convergence}
\end{figure}

In this case, the velocity field as well as the void fraction are unsteady and the velocity field is not divergence free. This leads to a nonzero mass source that should be accounted for. A very fine mesh size ($\Delta x = \Delta y = 0.015625 m$) is used in this case so that spatial errors are negligible. This is important for the convergence analysis in time. 
The velocity, pressure and void fraction fields at $t = 1s$ are shown in Fig. \ref{fig:case3} for a timestep of $0.001 s$ as to clearly and accurately represent the pressure field since the analytical solution is zero. The convergence plot is shown in Fig. \ref{fig:case3_convergence}.

\subsection*{Discussion}
The plots presented in Figs. \ref{fig:case1_convergence},  \ref{fig:case1_convergenceA}, and \ref{fig:case2_convergence} show that the velocity and pressure converge spatially to the correct order for all finite element order investigated. This demonstrates the high-order capabilities of the proposed stabilized scheme for the VANS equations. Due to longer simulation times, no scheme higher than Q3-Q3 is investigated. However, a similar behavior is expected for higher-order schemes. Since Lagrange polynomials are used as an interpolation basis, we do not expect the present code to be well conditioned for very high-orders ($>6$). In Fig. \ref{fig:case3_convergence} (a), we show that time convergence of our scheme is first-order for the different spatial schemes used. This is because time is discretized using the first order Backward Differentiation Formulation (BDF1). Similarly, in Fig. \ref{fig:case3_convergence} (b), we show second order convergence in time when using BDF2 as this represents a second order time discretization scheme. Thus, our code preserves the accuracy of the temporal discretization used.  Higher order temporal schemes are also compatible with the code. Spatial convergence for the transient case is not shown since exact convergence is not attained for the finer meshes. We attribute this discrepancy to the fact that a very small time step should be used in the study of spatial convergence of the transient case to reduce temporal errors. For the MMS solution, this leads to the presence of a very small time step. Using MMS, we determined that our solver is globally mass conservative. The maximum global losses is of the same order of magnitude as the tolerance of the non-linear solver. For steady flows, relative mass loss of the order of $10^{-17}$ could be obtained, which demonstrates the good mass conservation capabilities of this scheme.
For the transient case, the global mass losses ranges between $10^{-6}$ and $10^{-8}$. This error is dependent on the time discretization scheme as well as on the tolerance of the non-linear solver. \\


\section{Simulation of Packed Bed for the Validation of the Stabilized VANS Scheme}
After verifying our model using the method of manufactured solutions, we validate our model through the simulation of a packed bed and the prediction of the pressure drop inside it.

\subsection{Simulation Setup}
In order to realize this packed bed study, we perform two different simulations. The first simulation involves the discrete element method (DEM) which fills the bed with particles in order to create a non-uniform packing. The second simulation involves the VANS equations where we simulate a flow of fluid across the packing. A flow through a packed bed example on the documentation page of Lethe clearly represents the simulation setup with the different parameter files. All parameters necessary for the simulations can be found in Table \ref{tab:geometry}.

\subsubsection{Geometry Description}
The simulation is carried out in a 3D lab-scale cylindrical packed bed using Q1-Q1, Q2-Q1 and Q2-Q2 elements. The geometry is generated using the subdivided\_cylinder function of Deal.II. It has a total length of 200 $mm$ and a radius of 25 $mm$. The bed is divided into 20 divisions across its half length. The width of the bed walls is neglected. A stopper is defined at a distance of 45 $mm$ below the centroid of the bed to allow the suspension of particles and prevent their fall to the bottom of the bed. We used a non-uniform grid containing 3200 cells to describe the system. This is done to ensure better accuracy while respecting that the particle size should be at least one third smaller than the cell size in CFD simulations \cite{PengZhengbiao2014Iovf}. 

\subsubsection{Discrete Element Method (DEM)}
We perform the DEM simulations using Lethe-DEM \cite{GolshanShahab2021L:Ao,golshan2022load}. For the purpose of validation, 200,000 particles having the same physical properties were considered for this simulation. All parameters for the DEM simulation are shown in Table \ref{tab:geometry}. The particles have a diameter of 1 $mm$ and a density of 2500 $kg/m^3$. The particles are inserted randomly above the stopper where they fall under the influence of gravitational forces. We use the velocity Verlet method as the integration scheme for the motion of the particles. We used non-linear viscoelastic particle-particle and particle-wall contact force models. This is because a non-linear viscoelastic model presents a higher accuracy in calculating the collision forces \cite{norouzi2016}. The particles are allowed to settle on the stopper until their velocity dissipates. The DEM simulation is performed in parallel on 8 processors for 2 $s$ with a time step of $2\times10^{-5} s$. We are left with a packed bed column with a height of 80 $mm$. 

\subsubsection{VANS-CFD Simulation}
The gas chosen for this simulation is air. We define the various parameters of the VANS-CFD simulation in Table \ref{tab:geometry}. The gas is introduced with a uniform velocity distribution at the inlet located at the bottom of the cylindrical bed. The simulation is carried out in parallel on 8 processors for 0.5 $s$ with a time step of $2\times10^{-3} s$. For all cases, we ensure that steady state was reached at the end time. For the calculation of the fluid governing equations, we introduce a Dirichlet boundary condition for the velocity at the inlet of the bed. A slip boundary condition is applied to the walls of the cylinder. As an initial condition, the velocity is set to zero everywhere except at the inlet. We run the simulation for different inlet velocities which ranged between 0 $m/s$ and 0.6 $m/s$ with an increment of 0.05 $m/s$. This ensures that we remain below the minimum fluidization velocity $(U_{mf})$. There exist various correlations to predict the minimum fluidization velocity. The choice is usually based on the particles' diameter and density. For a more detailed review of the different available correlations, we refer the reader to the article by Anantharaman et al. \cite{ANANTHARAMAN2018454}. The minimum fluidization velocity $(U_{mf})$ was calculated based on the Wen-Yu correlation for Geldart D particles and is defined as \cite{WenC.Y1966Agmf}:
\begin{align}
    U_{mf}= \Bigg[\Big(33.7^2 + 0.0408 Ar\Big)^{0.5} - 33.7\Bigg] \frac{\mu_f}{\rho_f d_p}
    \label{eq::umf2}
\end{align}
where $Ar$ is the Archimedes number and is defined as:
\begin{align}
    Ar = \frac{g \rho_f (\rho_p - \rho_f) d_p^3}{\mu_f^2}
\end{align}
and $\rho_p$ is the particles' density, $g$ is the gravitational acceleration, and $\phi$ is the sphericity of the particle. $\phi=1$ for spherical particles. In our case, $U_{mf} = 0.718$ $m/s$. For the purpose of determining $U_{mf}$, the void fraction $\epsilon_f$ of the bed is calculated as:
\begin{align}
  \epsilon_f = \frac{(\pi r_b^2)H_b - n_p(\frac{4}{3}\pi r_p^3)}{(\pi r_b^2)H_b }  
\end{align}
where $n_p$ and $r_p$ are the particles' number and radius respectively, and $r_b$ and $H_b$ are the bed's radius and height respectively. From the DEM simulation, the bed height was found to be 89.86 $mm$.

\subsection{Results and Discussion}
\begin{figure}[ht]
         \centering
         \includegraphics[width=0.7\textwidth]{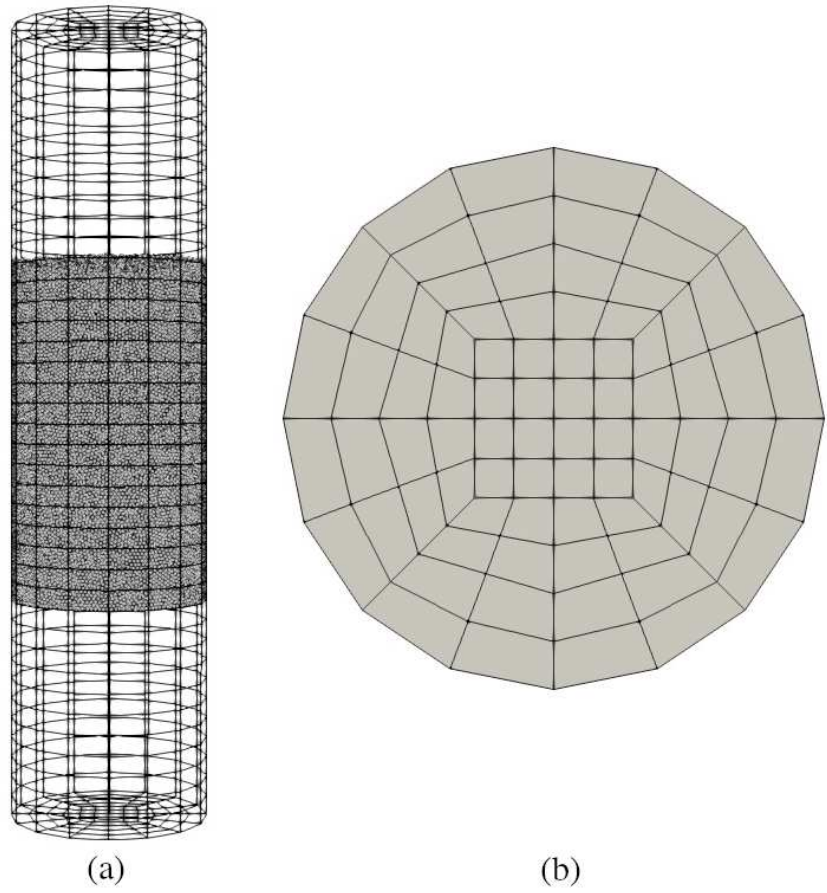}
        \caption{xy plane particle packing and mesh (a) and yz plane mesh (b) of packed bed}
        \label{fig:pb_dem}
\end{figure}

The particle packing result of the DEM simulation along with the meshed packed bed are shown in Fig. \ref{fig:pb_dem}. The particles take around 1.5 $s$ to settle on the stopper. However, the simulation was run for 2 $s$ in order to ensure that all particles become at rest. It is important that all particles possess a negligible velocity as the drag force is calculated using the relative velocity. If the particles had a non negligible velocity, the calculated drag would become incorrect resulting in an incorrect prediction of the pressure drop along the bed.

\begin{figure}[ht]
     \centering
         \includegraphics[width=0.7\textwidth]{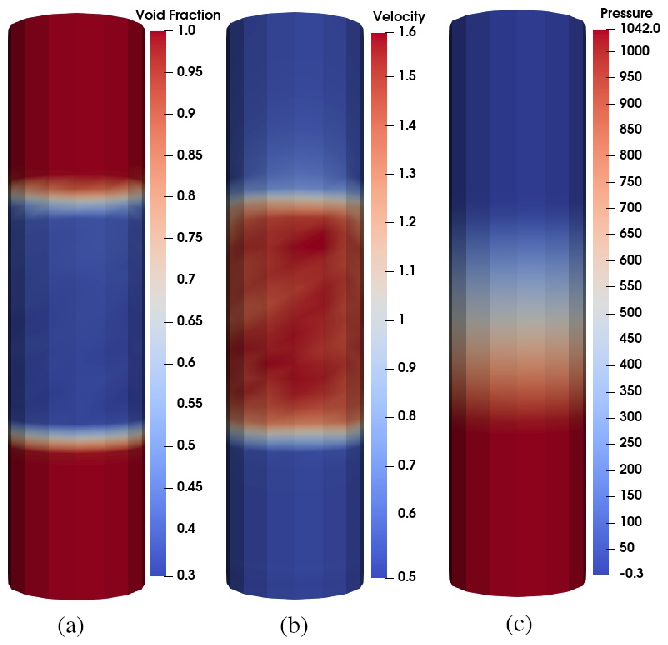}
         \caption{void fraction distribution}
        \caption{Solution of gas flow through packed bed using Q1-Q1 elements for (a) void fraction, (b) velocity and (c) pressure.}
        \label{fig:pb}
\end{figure}

The velocity, pressure, and void fraction distributions are presented in Fig. \ref{fig:pb}. The results shown are based on an inlet velocity $u_{in} = 0.5$ $m/s$, and the Rong et al. drag model \cite{RongL.W2013Lsof}. The void fraction shown in Fig. \ref{fig:pb} (a) appears to be spatially continuous. This is due to the smoothed $\mathcal{L}^2$ projection which overcomes the discontinuity of the Particle Centroid Method. This discontinuity is controlled using the smoothing length $L$. In this simulation, $L^2$ is taken to be equal to $5$ $d_p^2$. For the velocity distribution shown in Fig. \ref{fig:pb} (b), we can see some oscillations inside the packing region. This is explained by the non-uniform spacing between particles, resulting in higher velocities in regions with smaller spacing and vice versa. In Fig. \ref{fig:pb} (c), the pressure drops instantly inside the packing as expected. This is mainly due to friction inside the packing.

In order to validate our model, the results of our simulations are compared with the predicted $\Delta p$ of the Ergun equation \cite{Ergun} as a function of Reynolds number. This equation is a correlation used to predict the pressure drop $\Delta p$ across a packed bed of spherical particles and is denoted as:
\begin{align}
    \Delta p = \frac{150(1-\epsilon_f)^2 \bm{u}_f \mu_f H_b}{\epsilon_f^3 d_p^2} + \frac{1.75(1-\epsilon_f)\rho_f \bm{u}_f H_b}{\epsilon_f^3 d_p}
    \label{eq::ergun}
\end{align}
The Reynolds number used is the bed's Reynolds number:
\begin{align}
    Re = \frac{\rho_f \bm{u}_f D}{\mu_f}
\end{align}
where $D$ is the bed's diameter.

\begin{figure}[ht]
        \centering
        \includegraphics[width=0.7\textwidth]{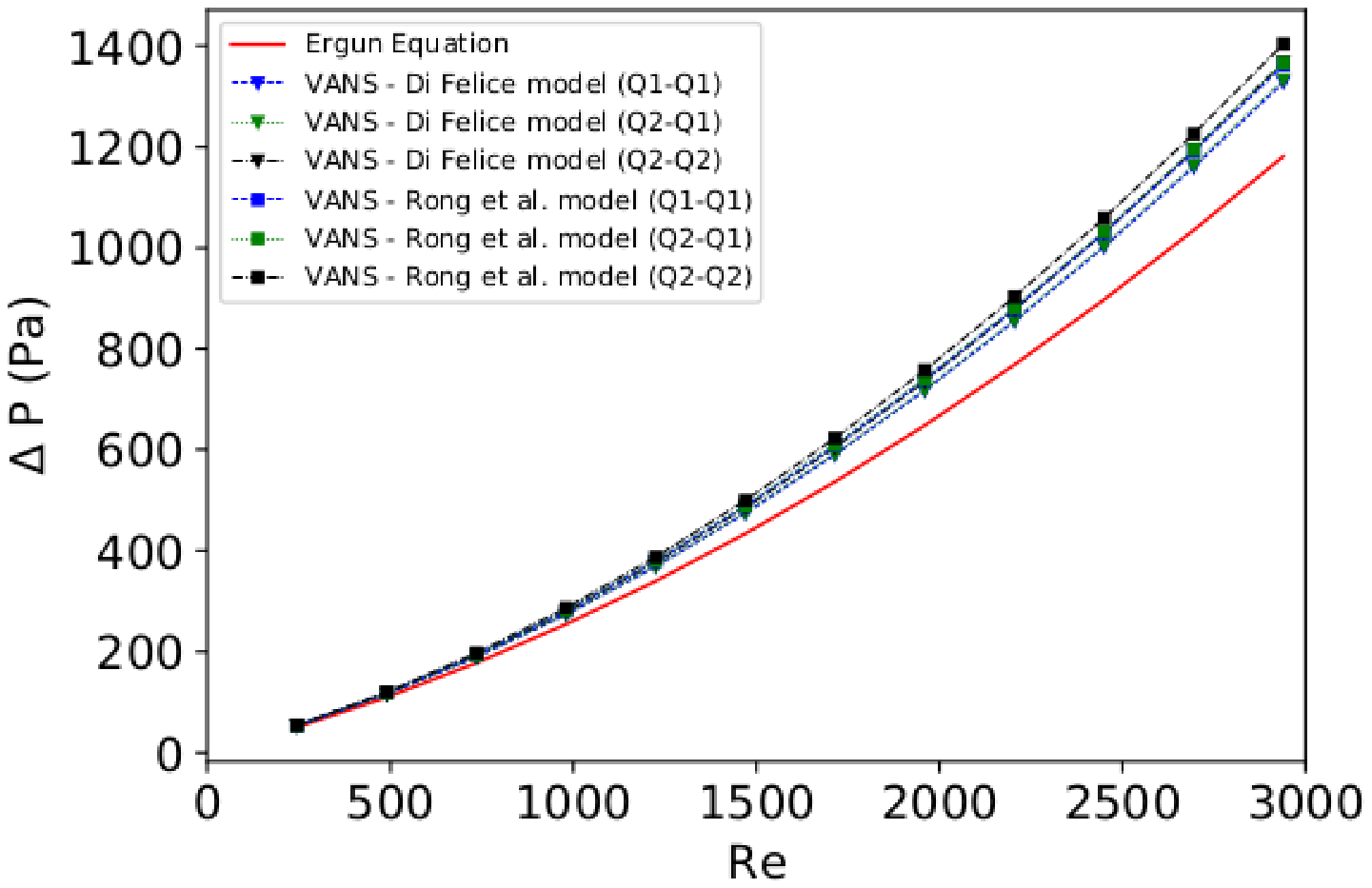}
        \caption{Comparison of model B pressure drop in packed bed as a function of inlet Reynolds number for different drag models and finite element orders.}
        \label{fig:delta_p}
\end{figure}
\begin{figure}[ht]
        \centering
        \includegraphics[width=0.7\textwidth]{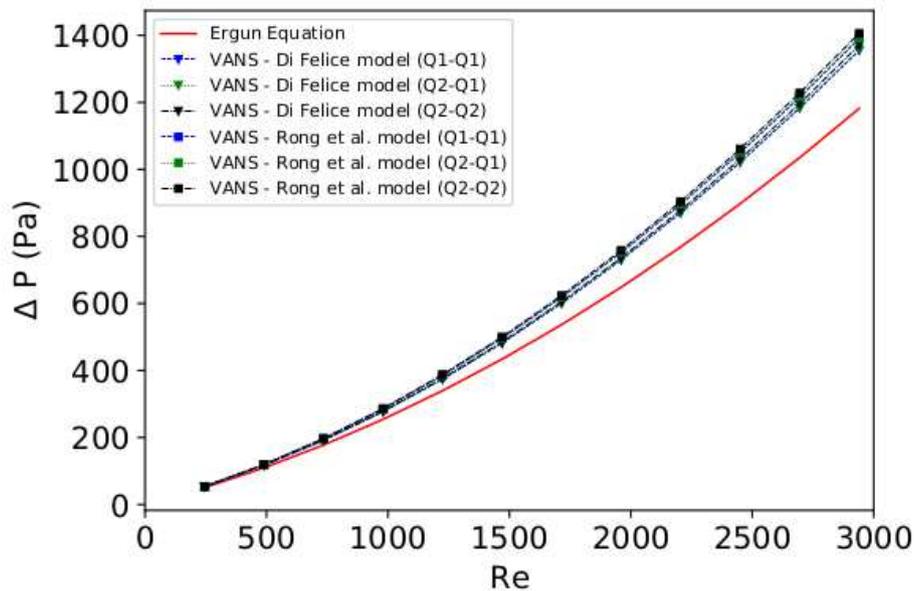}
        \caption{Comparison of model A pressure drop in packed bed as a function of inlet Reynolds number for different drag models and finite element orders.}
        \label{fig:delta_pA}
\end{figure}

For the Ergun equation, $H_b$ is calculated in the same way as for the minimum fluidization velocity ($U_{mf}$). Similarly, in this case, $H_b = 89.86 mm$. Figs. \ref{fig:delta_p} and  \ref{fig:delta_pA} compare the pressure drop curves obtained from the simulations using Di Felice \cite{DiFeliceR1994Tvff} and Rong et al \cite{RongL.W2013Lsof} models for Q1-Q1, Q2-Q1 and Q2-Q2 elements using models B and A respectively with the predicted pressure drop obtained from the Ergun equation.

Using model B of the VANS equations, the two drag models show a good correlation with Eq. (\ref{eq::ergun}). They both slightly over predict the pressure drop compared to the Ergun equation for higher $Re$. In both cases, the Rong et al. drag model shows more over prediction compared to the Difelice drag model. 
The $\beta$ coefficient in the Rong et al model takes into account not only the particle's Reynolds number ($Re_p$) but also the bed's void fraction ($\epsilon_f$). This can explain the over prediction of this model at higher $Re$ values where the effects of the void fraction are amplified by the Reynold's number. 
Moreover, for the same drag model, using higher order elements leads to a higher over prediction of the pressure drop at higher $Re$ values and does not change the predicted results at lower $Re$ values. 
Using model A of the VANS equations, both drag models give similar results for the different element orders. As the order of the elements increased, the variation of the pressure drop remained negligible compared to the values obtained from Model B. However, the values of the pressure drop obtained using Model A were in general slightly greater than those obtained using Model B for the same element order.
These results prove that the drag models' implementation in the VANS solver is correct. The local calculation of the drag in each element along with the volumetric integral of its value over all the elements successfully predicts the pressure drop in the packed bed. It also demonstrates that higher order elements are better suited for efficiently solving problems with a requirement for high accuracy, but for which the minimum mesh size is limited.

\begin{figure}[ht]
        \centering
        \includegraphics[width=0.7\textwidth]{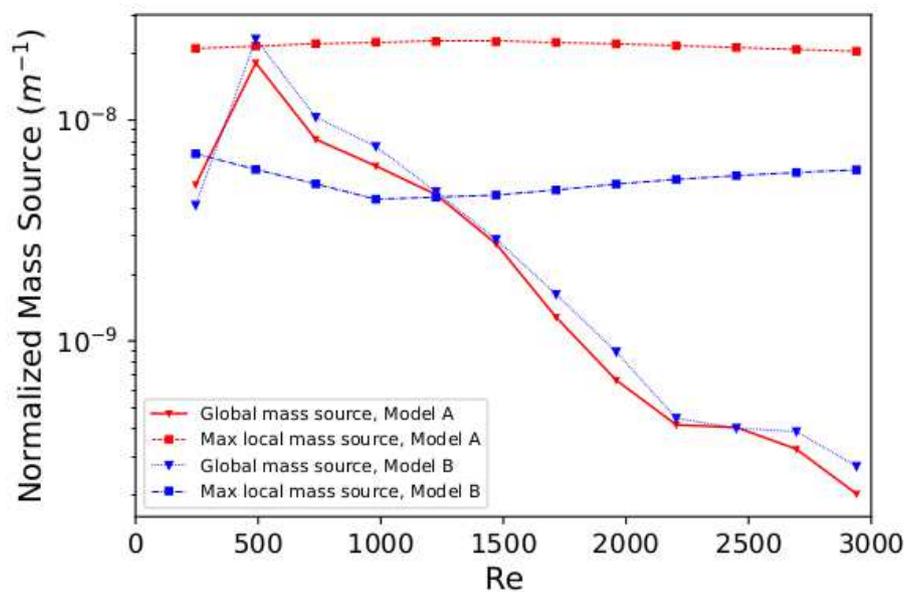}
        \caption{Comparison of local and global mass conservation of models A and B of the VANS equations using Rong drag model and Q1-Q1 elements.}
        \label{fig:mass_source}
\end{figure}

For the same packed bed simulation without bounding the void fraction, we determined the global mass source term and the maximum local mass source term at every inlet velocity for the Rong et al. \cite{RongL.W2013Lsof} drag model using both models A and B of the VANS equations. The maximum local mass source term was determined by integrating the continuity equation over each element and choosing the value of the element with the maximum absolute value. The global mass source term was determined by summing up all element-wise integration of the continuity equation. The mass source terms are shown in absolute value scaled by the inlet velocity as in Fig. \ref{fig:mass_source}.

Both models A and B show the same trend as well as similar values of the global mass source for different inlet velocities. Therefore, both our models are globally mass conservative achieving mass source term errors in the magnitude of $10^{-10}$, which is the order of magnitude of the tolerance of the non-linear solver. Moreover, the two models show similar behavior for the maximum local mass conservation. As the inlet velocity increases, the local mass source for both models remains relatively constant in the magnitude of $10^{-9}$ for model B and $10^{-8}$ for model A. Even though Model B achieves better local mass conservation, the magnitude of the local error obtained for both models is relatively small. Therefore, it is safe to consider that without bounding the void fraction, both models are mass conservative. In this case, bounding the void fraction does not add any benefits to the simulation as the projected void fraction had a maximum value of 1.00001 and thus it didn't affect any aspect of the simulation.  Some of the simulations were repeated while bounding the void fraction. The obtained results show diminished mass conservation in the magnitude of $10^{-6}$. As we bound the void fraction, we slightly change the values of the void fraction at the extremities of the packing and as a consequence alter the mass preservative properties of the $\mathcal{L}^2$ projection. This explains the increase in the magnitude of the mass source term error.

\section{Conclusion}
This work presents a verified and validated stabilized finite element approach for the modeling of both forms A and B of the volume-averaged Navier Stokes equations which are used to model multi-phase flows. We were able to develop a parallel high-order finite element solver for these equations within Lethe, thus enhancing accuracy without requiring finer meshes. This is important as it helps circumvent one of the solid-fluid flow limitations concerning mesh size. We verified our solver using the method of manufactured solution and showed that our model is consistent with the underlying discretization and globally mass conservative. The validation of the model through the prediction of the pressure drop in a packed bed simulation and its comparison with that of the Ergun correlation allowed us to confirm the validity of the drag models' implementation. In the presence of a discontinuous void fraction field, we controlled the stability of the simulation through our choice of the smoothing length $L$. This is often taken to be two or three times the particles' diameter. 
Finally, the methodology introduced here serves as a first step in developing a fully parallel high-order CFD-DEM solver by coupling the VANS solver presented in this paper with Lethe-DEM \cite{GolshanShahab2021L:Ao}. The rigorous verification and validation strategy deployed in this work can also be followed by other researchers who are developing solvers for the VANS equations, without regard to the numerical method used.


\section{Acknowledgements}
The authors would like to thank the deal.II community for their support. Bruno Blais would like to acknowledge the financial support from the Natural Sciences and Engineering Research Council of Canada (NSERC) through the RGPIN-2020-04510 Discovery Grant. The authors would also like to acknowledge technical support and computing time provided by Compute Canada and Calcul Québec.

\bibliography{wileyNJD-AMA}%

\newpage

\begin{table}[ht]
 \captionof{table}{Physical and numerical parameters for the packed bed simulation}
 \label{tab:geometry}
\centering
\resizebox{\textwidth}{!}{%
 \begin{tabular}{l l l l} 
 \hline
 \textit{Simulation control}\\
 DEM end time ($s$) & 2 & VANS-CFD end time ($s$) & 0.5\\
 DEM time step ($s$) & $2\times10^{-5}$ & VANS-CFD time step ($s$) & $2\times10^{-3}$ \\
 \\
 \textit{Geometry}\\
 Bed Height ($mm$) & 200 &  Bed Radius ($mm$) & 25 \\
 Wall Thickness ($mm$) & 0\\
 \\
 \textit{Particles}\\
 Number & $2\times10^{5}$ & Diameter(m) & $1\times10^{-3}$\\
 Density($kg/m^3$) & 2500 & Young Modulus ($N/m^2$) &  $1\times10^6$\\
 Particle-particle poisson ratio & 0.3 & Particle-wall poisson ratio & 0.3 \\
 Particle-particle restitution coefficient & 0.2 & Particle-wall restitution coefficient & 0.2 \\
 Particle-particle friction coefficient & 0.1 & Particle-wall friction coefficient & 0.1 \\
 Particle-particle rolling friction & 0.2 & Particle-wall rolling friction & 0.3 \\
 \\
 \textit{Gas phase}\\
 Viscosity ($Pa\cdot s$) & $1\times10^{-5}$ & Density ($kg/m^3$) & 1 \\
 Inlet Velocity ($m/s$) & [0.05-0.6] & Void fraction smoothing factor $L^2$ & $5\times10^{-6}$\\
 \\
 \textit{Linear Solver}\\
 Method & GMRES & Max iterations & 5000\\
 Minimum residual & $1\times10^{-11}$ & Relative residual & $10
 \times10^{-3}$\\
 ILU preconditioner fill & 1 & ILU preconditioner absolute tolerance & $1\times10^{-14}$\\
 ILU preconditioner relative tolerance & 1\\
 \\
 \textit{Non-linear solver}\\
 Tolerance & $1\times10^{-9}$ & Max iterations & 10\\
 \hline
 \end{tabular}}
\end{table}

\newpage

\appendix
\section{Grad-Div Stabilization}
\begin{figure}[ht]
     \centering
        \includegraphics[width=0.7\textwidth]{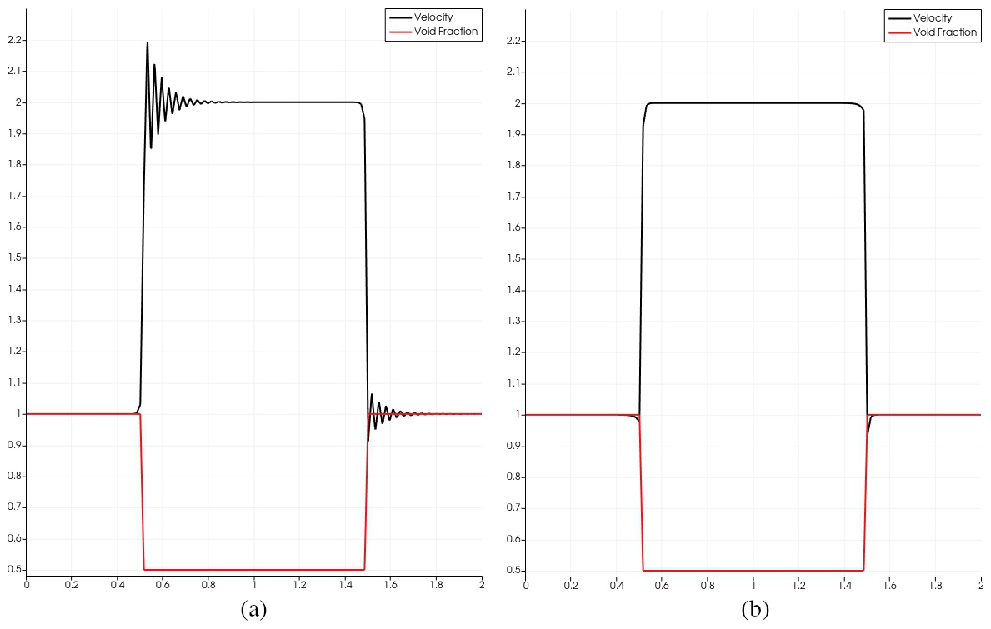}
        \caption{Comparison of simulation results of a fluid flow through a section with a sharp void fraction gradient with SUPG stabilization only (a) and with SUPG and grad-div stabilization (b).}
        \label{fig:appendix}
\end{figure}

In this appendix, we demonstrate the effect of proper grad-div stabilization on the stability of the VANS solver in the presence of a sharp void fraction gradient. A step function for the void fraction is defined as:
\begin{align}
    \epsilon_f = 
    \begin{cases}
      1 & \text{if $x < -0.5$ or $x > 0.5$}\\
      0.5 & \text{if $ -0.5 \leq x \leq 0.5$}
      \end{cases}
\end{align}

and an inlet velocity $\bm{u} = 1$ is defined in a square domain $[-1,1]^2$. The simulation is performed with and without the grad-div stabilization. The results obtained are shown in Fig. \ref{fig:appendix}

The simulation with only SUPG and PSPG stabilization is not enough to dampen the spurious oscillations in the velocity field arising from the sharp void fraction gradients at the location of the discontinuities. When grad-div stabilization is implemented, the oscillations disappear which makes the simulation much more stable. For our solver, using this grad-div stabilization method alleviates the need of a shock capturing scheme to diffuse the spurious oscillations.

\end{document}